\documentclass[conference]{IEEEtran}
\usepackage{amssymb}
\usepackage{amsmath}
\usepackage{graphicx}
\usepackage{algo}

\ifCLASSINFOpdf
\else
\fi
\begin{document}
%
\title{Primal Dual Affine Scaling on GPUs}
\author{\IEEEauthorblockN{Nithish Divakar}
\IEEEauthorblockA{Supercomputer Education and Research Center\\
Indian Institute of Technology, Bangalore\\
nithish.divakar@ssl.serc.iisc.in}
}
%
\maketitle
\begin{abstract}
Here we present an implementation of Primal-Dual Affine scaling method to solve linear optimisation problem
on GPU based systems. Strategies to convert the system generated by complementary slackness theorem into a symmetric system are given. A new CUDA friendly technique to solve the resulting symmetric positive definite subsystem is also developed. Various strategies to reduce the memory transfer and storage requirements were also explored.

\textit{Keywords|}Primal Dual Affine scaling, Linear programming, Interior point algorithms, Woodbury formula
\end{abstract}


%
\section{Introduction}
	A Linear Programming Problem  concerns with maximising or minimising an objective function subjected to certain constraints. Canonical form of a linear programming problem can be written as follows.
\begin{align}
	\label{eq:primal}
	minimize~&c^Tx \nonumber\\
	&Ax=b\\
	&x\geq0 \nonumber
\end{align}
The form given in (\ref{eq:primal}) is called primal problem and its corresponding dual problem is expressed as 
\begin{align}
	\label{eq:dual}
	maximize~&b^Ty \nonumber\\
	&A^Ty+s=c\\
	&s\geq0 \nonumber
\end{align}
\begin{align*}
	A\in & \mathbb{R}^{m\times n}\\
	x\in & \mathbb{R}^n,~y\in  \mathbb{R}^m \text{ is set of basic variables}\\
	c\in & \mathbb{R}^n,~b\in  \mathbb{R}^m\\
	s\in& \mathbb{R}^n\text{ is the slack variable}
\end{align*}
	
In equation (\ref{eq:primal}) and (\ref{eq:dual}), A is assumed to have full rank and $m < n$.
	A solution to the primal problem is a vector $x^* \in \mathbb{R}^n$ such that $c^Tx^*$ is minimised
and the constraints $Ax^* = b$, $x^*\geq 0$ are each satisfied.
\section{Solving Linear Programming Problems} 
	The constraints of a linear optimisation problems forms an n dimensional polyotope and the objective function is an n-1 diminutional plane in the same space. The minimum of objective function is has its minimum at one of the vertices of this constraint polyotope. This fact was proved by Chong \cite{Chong:2008fk}. most of the linear programming methods are based on searching for this vertex from an initial point in the polyotope.

According to the way in which the searching progresses, the methods to solve linear programming problems can be broadly divided into Boundary methods and Interior point methods. Boundary methods always stays in the boundary of the polyotope and search for the optimum vertex where the objective function is minimum. 

Simplex method, the most common boundary point methods start with an initial feasible solution and move towards adjacent vertices, constantly reducing value of objective function in the process till it cannot be reduced further. A vivid account of how the method works and the theory behind it can be found in Chong \cite{Chong:2008fk} and Nazareth \cite{Nazareth:2004}

Interior point methods start by selecting a solution inside the polyotope and iteratively advance the point towards the optimum node. The fastest methods pick initial point in the neighbourhood of center of the polyotope and follow a path central from all nodes towards the optimum nodes. Detail account of this can be found in Saigal \cite{Saigal:1995}

\begin{center}
	\includegraphics[width=6cm]{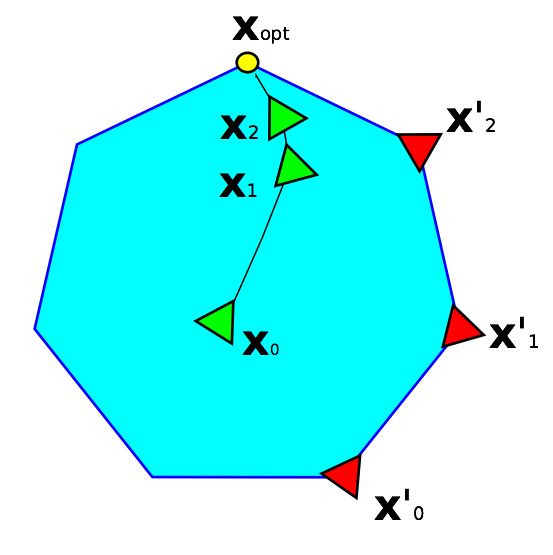}\\
	Polytope of a two-dimensional feasible region. The $x_i$ points are determined by an interior point method, while the x$_i^\prime$ ones by a simplex method.
\end{center}

\section{The Algorithm}
By complementary slackness theorem in Dantzig and Orden \cite{Dantzig:1953} , the following system of non-linear is generated for the primal and Dual problems described in equations (\ref{eq:primal}) and (\ref{eq:dual}).
$$
\begin{tabular}{ccccc}
Ax&&&=&b\\
&A$^T$y+&s&=&c\\
&&Xs&=&0
\end{tabular}
$$
On applying newtons method to this system we can find the optimum direction to which the next step is to be taken. This is given by
\begin{equation}\label{eq:bigsystem}
	\begin{bmatrix}
	A&O&O\\
	O&A^T&I\\
	S&O&X
	\end{bmatrix}\begin{bmatrix}
	\Delta x\\
	\Delta y\\
	\Delta s
	\end{bmatrix}=\begin{bmatrix}
	\vec{0}\\
	\vec{0}\\
	-Xs
	\end{bmatrix}
\end{equation}

Where S and X are diagonal matrices with diagonal entries of s and x respectively.

The primal dual affine scaling algorithm from Saigal \cite{Saigal:1994} can be outlined as follows:\\
\begin{center}
\begin{minipage}{10cm}
\begin{algorithm}{Primal-Dual Affine Scaling}
	\Step{Let $(x^0,y^0,s^0)$ be an initial interior solution}
	\Step{0 $< \alpha <$ 1}
	\Step{k:=0}
	\While{TRUE}
		\Step{Solve system (\ref{eq:bigsystem}) to obtain optimal}
		\Step{direction vectors ($\Delta x^k,\Delta y^k,\Delta s^k$)}
		\If{$\Delta x^k=0,\Delta y^k=0$ and $\Delta s^k=0$}
			\Step{$(x^k,y^k,s^k)$ is Optimum. STOP}
		\EndIf
		\Step{$\psi_k=max\left\{\phi(-X_k^{-}1\Delta x^k),\phi(-S_k^{-1}\Delta s^k)\right\}$}
		\Comment{$\phi$(u)=max\{u$_j$\} for any vector u}
		\If{$\psi_k=1$}
			\Step{$\alpha:=1$}
		\EndIf
		\Step{$x^{k+1}=x^k+\frac{\alpha}{\psi_k}\Delta x^k$}
		\Step{$y^{k+1}=y^k+\frac{\alpha}{\psi_k}\Delta y^k$}
		\Step{$s^{k+1}=s^k+\frac{\alpha}{\psi_k}\Delta s^k$}
		\If{$\psi_k=1$}
			\Step{$(x^k,y^k,s^k)$ is Optimum .STOP}
		\EndIf
		\Step{k=k+1}
	\EndWhile
\end{algorithm}
\end{minipage}
\end{center}
The most important part of this algorithm is in solving the system (\ref{eq:bigsystem}). As we can see from the algorithm itself, this system has to be solved in every iteration. This create a heavy overhead as all other steps in the algorithm are trivial. This work mainly concentrates on solving the system (\ref{eq:bigsystem}).
\section{Previous Work}
	The area of Linear optimisation is well developed in and has a strong theory base. Some of the earlier works in this field have made the subject evolve and reach its current state. Theoretically the fastest known interior point algorithm to date is due to Vaidya \cite{Vaidya:1990} which requires $O(((m + n)n^2+ (m + n)^{1.5}n)L)$ arithmetic operations where m is the number of constraints, and n is the number of variables. The algorithm has been reported to work in O(L) precision. Also the work of Mehora \cite{Mehora:1992} on implementation of primal-dual methods are encouraging. In his work, the primal-dual method have been approximated by second order terms to provide faster convergence and use of potential functions to find direction and stay close to central path along iteration paved the way for development of central path based algorithms. Goldberg et.al \cite{Goldberg:1989} describes the use of interior point methods to find solution of bipartite matching and minimum cost related problems. It shows that Linear optimisation can be used to solve various other problems from different domains.

Recent developments in GPU technology has allowed the growth of Linear programming into this domain as well. Spampinato and Elster \cite{Elster:2009} implemented a variant of simplex method on NVIDIA GPUs using CUDA programming model. A speedup of 2x-2.5x was reported. Work of Smith et.al \cite{Smith:2011} on Linear optimisation problems with unknown constraint matrix also opened up a new  area of research. They have used accelerated sparse matrix-vector product computation to accelerate convergence rate. 
\section{The Proposed Strategy}
	The major bottleneck of the Primal-Dual Affine Scaling algorithm is in solving the system of equations (\ref{eq:bigsystem}).  Practical Optimisation problems consist of thousands of constraints and millions of variables. As such the matrix A is quiet large and system (\ref{eq:bigsystem}) is even larger. So we propose the following strategy.

Simple Row column interchange will transform system (\ref{eq:bigsystem}) into following symmetric system 
\begin{equation}\label{eq:bigsystem2}
	\begin{bmatrix}
		O&A^T&I\\
		A&O&O\\
		I&O&D
	\end{bmatrix}
	\begin{bmatrix}
		\Delta x\\
		\Delta y\\
		\Delta s
	\end{bmatrix}=
	\begin{bmatrix}
		\vec{0}\\
		\vec{0}\\
		-x
	\end{bmatrix}
\end{equation}
Where $D=S^{-1}X$ is an $n\times n$ diagonal matrix.
This gives the components as 
\begin{align*}
	\Delta x=&D[-\Delta s]-I\vec{x}\\
	(ADA^T)\Delta y=&A\vec{x}\\
	\Delta s=&-A^T\Delta y
\end{align*}
Our main concern is in solving the system $(ADA^T)\Delta y=A\vec{x}$.
\subsection{Solving $(ADA^T)x=b$}
These type of systems are usually solved using cholesky factorisation as $(ADA^T)$ is symmetric and positive definite. But this strategy is not optimum for architecture like CUDA as the decomposed matrix has a triangular shape easing to lots of overhead to coordinate and idle threads. More over, the triangular system obtained at the end gives no scope of parallelism at all if the system is dense.

A better optimisation strategy involves the observation that only the D matrix in $ADA^T$ changes in every iteration. Thus, if we can develop a strategy to split the matrix into a constant part and a variable part, then lots of computation can be done before iteration begins. Such a strategy is obtained from Woodbury formula.
$$(\mathbb{A}+UV^T)^{-1}=\mathbb{A}^{-1}+G\mathbb{D}G^T$$
If $\mathbb{A}=AA^T$, $U=A$ and $V=A(D-I)$ then $(\mathbb{A}+UV^T)^{-1}=(ADA^T)^{-1}$. When the change in matrix is considered as series of rank one updates, the following formula can be derived for directly finding solution of the system.

\begin{align*}
	x_l=&x_{l-1}-\left[\frac{v_l^Tx_{l-1}}{1+v_l^Ty_{l-1,l}}\right]y_{l-1,l}\\
	y_{l,k}=&y_{l-1,k}-\left[\frac{v_l^Ty_{l-1,k}}{1+v_l^Ty_{l-1,l}}\right]y_{l-1,l}
\end{align*}
\begin{flushright}
$k=l,l+1,\ldots,n$
\end{flushright}
Where $Y=\left[y_{0,1}|y_{0,2}|\cdots|y_{0,n}\right]=(AA^T)^{-1}A$ and $x_0=Yb$. $x_n$ is the solution of the system $(ADA^T)x=b$. This formulation is due to Egidi and Maponi\cite{maponi2006}.
\subsection{CUDA Related Optimisations}
\subsubsection{Data Storage}
The x vector can be appended as a column with the Y matrix so that the operations can be done uniformly. Another requirement was to store the inner products. So the following distribution is adopted.
\begin{center}
	\begin{tabular}{|ccccc|c|}
		\hline
		&&&&&\\
		&&Y&&&x\\
		&&&&&\\
		\hline
		\multicolumn{6}{|c|}{Inner Products}\\
		\hline
	\end{tabular}
\end{center}
V matrix entries are only accessed once and hence we can refrain from storing it by calculating the required entry directly as follows:
$$V_{i,j}=A_{i,j}(D_{i,i}-1)$$
\subsubsection{Partitioning for thread allocation}
Since there is no dependencies between columns, we will follow a partitioning strategy where each column is assigned to one block and inside the block, one thread takes care of only one element. When no of elements in column is more than maximum limit of threads, a 1D cyclic pattern of thread to element assignment is followed inside a block. Figure 1 describes the assignment of blocks to column.
\begin{figure}[htbp]
\begin{center}
	\includegraphics[width=4cm]{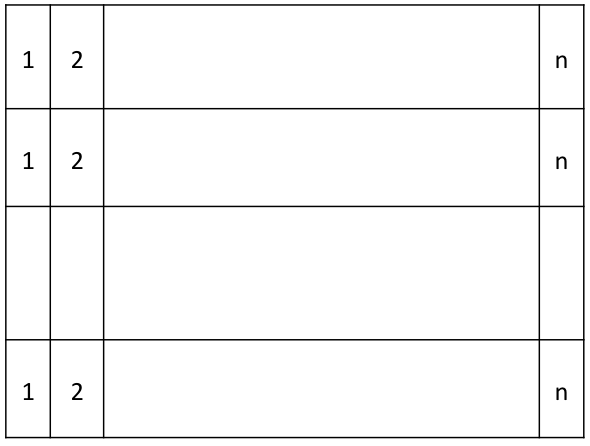}
\end{center}
	\label{fig_arrangement}
	\caption{Assignment of blocks to matrix}
\end{figure}
\subsubsection{Memory access pattern}
As the adjacent threads is assigned to adjacent elements of a column, the access pattern is strided by the matrix size in x direction. So to regularise memory access pattern, the whole matrix is stored in its transposed format. This transformation gives coalesced memory access  throughout out the  algorithm which is highly desirable for an architecture like CUDA.
\subsubsection{Inner product computation}
This  operation can be optimised for execution in CUDA  by following strategies:
\begin{itemize}
	\item The data is first read into shared memory and reduction is done only from shared memory
	\item Using the \textit{Add During Load} technique, idle threads can be eliminated. Moreover, total number of threads and required amount of shared memory can be reduced into half. This allows the algorithm to work with larger data set.
	\item Last 5 iteration can be unrolled and explicit synchronisation can be avoided because the threads in a warp has implicit synchronisation between them.
\end{itemize}
\subsubsection{Kernel Calls}
Multiple kernel calls are used to produce inter block synchronisation. To avoid kernel call overhead, the concept of constant memory is used. The argument values of all the kernels is copied to the constant memory residing in the device. Kernels directly access the arguments from constant memory and hence passing argument from host to device can be eliminated. Moreover parameter reading from constant memory is fast operation as all threads would simultaneously accessing same data which can be satisfied with just one read inside a warp. Moreover internal caching mechanism removes the over head of multiple reads. Hence this strategy keeps kernel call overhead to a minimum
\section{Experimental Results}
	The code has been executed in Tesla cluster  with following technical details
\begin{itemize}
	\item Each node of the cluster consists of Four AMD Quad-Core Opteron 8378 processors with 2.4Ghz clock speed.
	\item Each of these compute nodes is also connected to a NVIDIA-Tesla S1070 GPGPU node
 	\item Each Tesla node is composed of 4 GPUs with each GPU made up of 240 processor cores.
	\item 64GB Main Memory
	\item Nvidia Tesla S1070 1U server with 4 GPU's operating at 1.296Ghz
\end{itemize}
\begin{figure}[h!]
	\begin{center}
		\includegraphics[width=9cm]{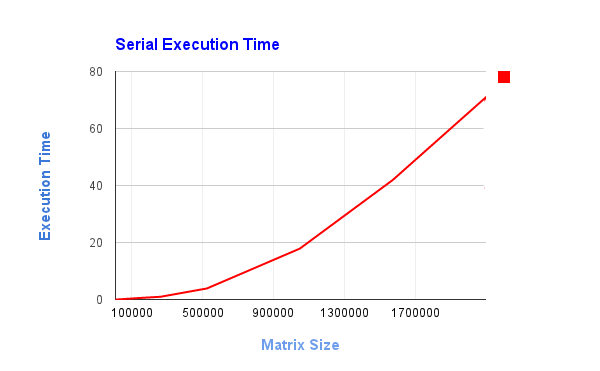}
	\end{center}
	\vspace{-0.4cm}
	\caption{Serial code execution time}
\end{figure}
\begin{figure}[h!]
	\begin{center}
	\includegraphics[width=9cm]{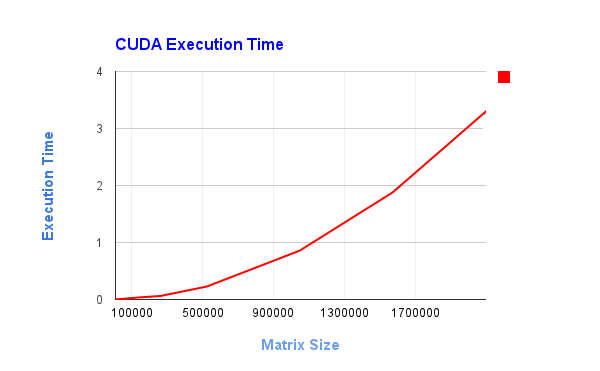}
\end{center}
	\vspace{-0.4cm}
	\caption{Execution time of the code implemented in CUDA}
\end{figure}
\begin{figure}[h!]
	\begin{center}
		\includegraphics[width=9cm]{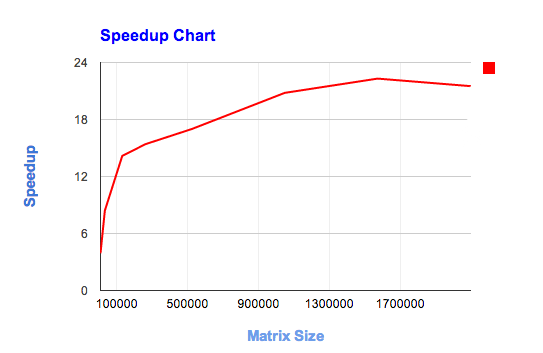}
	\end{center}
	\vspace{-0.4cm}
	\caption{This chart shows the speedup of the parallel code execution compared to the serial version. Highest observed speedup was 22x}
\end{figure}
\begin{figure}[h!]
	\begin{center}
		\includegraphics[width=7cm]{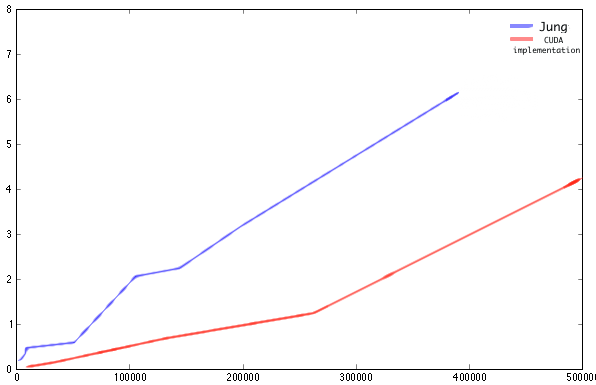}
	\end{center}
	\vspace{-0.4cm}
	\caption{Comparison of the Execution time of the whole algorithm with the execution time given in Jung et. al \cite{jin2008} . A speedup of 2x is seen from this comparison.}
\end{figure}
\begin{figure}[h]
	\begin{center}
		\includegraphics[width=7cm]{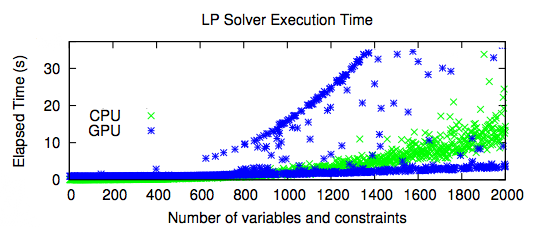}
	\end{center}
	\begin{center}
		\includegraphics[width=7cm]{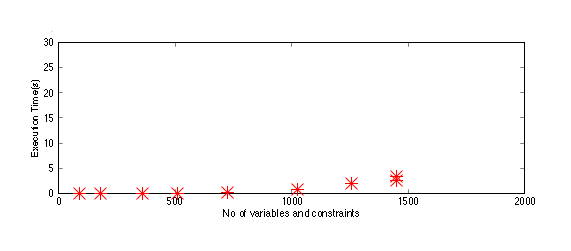}
	\end{center}
	\vspace{-0.4cm}
	\caption{Comparison of the Execution time between the method described in this work and the work done by Spampinato and Elster \cite{Elster:2009} on implementing simplex method on GPU}
\end{figure}
\section{Future Work}
	\begin{enumerate}
	\item The above formulation of the problem is for dense systems. This can be extended towards sparse systems. Maintaining the same degree of coalesced memory access will be challenge there.
	\item The application of Woodbury formula proved to be a good alternative for solving linear system for this kind of problems. The derived formula proved to be an optimum CUDA friendly algorithm for solving system of equations. This can be pursued further to be extended towards sparse systems. Maintaining sparsity of the system will be a challenge there.
	\item We have only tested this work in CUDA architecture with compute capability 1.0. This work needs to be tested with higher capabilities devising new techniques to use architectural features available there.
	\item Much more work has to be done in finding an efficient technique to reduce kernel call overhead. The technique describe here in this research effort proved to be fruitful one. Other options exist which can be pursued further.
	\item Although this formulation is successful in reducing computation time by totally eliminating the need for memory transfer in-between iterations, the area of hybrid execution can be explored if asynchronous memory transfers can be used.
\end{enumerate}

\bibliographystyle{plain}
\bibliography{bibliography}

\begin{thebibliography}{10}

\bibitem{Chong:2008fk}
Edwin K.~P. Chong and Stanislaw~H. Zak.
\newblock {\em An Introduction To Optimization}.
\newblock Jhon Wiley and Sons, 2008.

\bibitem{Dantzig:1953}
G.~B. Dantzig and A.~Orden.
\newblock Duality theorems.
\newblock Technical Report 1265, The RAND Corporation,Santa Monica, 1953.

\bibitem{Goldberg:1989}
A.~V. Goldberg, S.~A. Plotkin, D.~B. Shmoys, and E.~Tardos.
\newblock Interior point methods in parallel computation.
\newblock Technical report, Department of Computer Science, Stanford
  University, 1989.

\bibitem{jin2008}
Jin~Hyuk Jung and Dianne~P. O'Leary.
\newblock Implementing an interior point method for linear programs on a
  cpu-gpu systems, 2008.

\bibitem{Mehora:1992}
Sanjay Mehora.
\newblock On implementation od primal-dual interior point method.
\newblock In {\em Siam J. Optimization}. Soziety for Industrial Applied
  Mathematics, 1992.

\bibitem{maponi2006}
P.~Maponi N.~Egidi.
\newblock A sherman--morrison approach to the solution of linear systems, 2006.

\bibitem{Nazareth:2004}
John~Lawrence Nazareth.
\newblock {\em An Optimization Primer}.
\newblock Springer, 2004.

\bibitem{Saigal:1994}
Romesh Saigal.
\newblock On primal dual affine scaling method.
\newblock In {\em Opsearch}, 1994.

\bibitem{Saigal:1995}
Romesh Saigal.
\newblock {\em Linear Programming: A Modern Integrated Analysis}.
\newblock Kluwer Academic Publishers, 1995.

\bibitem{Smith:2011}
E.~Smith, J.~Gondzio, and J.~A.~J. Hall.
\newblock Gpu acceleration of the matrix-free interior point method.
\newblock Technical report, School of Mathematics and Maxwell Institute for
  Mathematical Sciences, The University of Edinburgh, 2011.

\bibitem{Elster:2009}
Daniele~G. Spampinato and Anne~C. Elster.
\newblock Linear optimization on modern gpus, 2009.

\bibitem{Vaidya:1990}
Pravin~M. Vaidya.
\newblock An algorithm for linear programming which requires o$(((m + n)n^2+ (m
  + n)^{1.5}n)l)$ arithmetic operations.
\newblock Technical report, AT\&T Bell Laboratories, 1990.

\end{thebibliography}
\end{document}